\documentclass[]{article}
\usepackage{amsmath}
\usepackage{amsfonts}
\usepackage{amssymb}
\title{On a Problem of Gromov about Generalizing Alexandrov-Fenchel Inequality}
\author{Yuri Burda}

\begin{document}

\maketitle

\begin{abstract}
In this note we give an answer to a question about mixed volumes asked by Gromov in  "Convex Sets and Kahler Manifolds". For reader's convenience we remind definitions and some of the properties of mixed volumes and mixed discriminants.

Dans cette note, nous donnions une r\'{e}ponse \`{a} une question sur les volumes mixtes pos\'{e}es par M. Gromov dans  "Convex Sets and Kahler Manifolds". Pour la commodit\'{e} du lecteur, nous rappelons les d\'{e}finitions et certaines des propri\'{e}t\'{e}s de volumes mixtes et discriminants mixte.

Mathematics Subject Classification: 52A39
\end{abstract}

\section{Mixed Volumes and Mixed Discriminants}

By a theorem of Minkowski the volume of a positive linear combination $\lambda_1 A_1+\ldots+\lambda_k A_k$ of $k$ convex bodies in $\mathbf{R}^n$ is a homogeneous polynomial of degree $n$ in $\lambda$'s: 
$$V( \lambda_1 A_1+\ldots +\lambda_k A_k)=\sum_{\substack{I\subset \mathbf{Z}_{+}^k \\ |I|=n}} {n \choose I} V_I \lambda^I$$

The coefficient $V_I$ for $I=\{i_1,\ldots,i_k\}$ is called the mixed volume $$V(\underbrace{A_1,\ldots,A_1}_{i_1},\ldots,\underbrace{A_k,\ldots,A_k}_{i_k})$$ of the bodies $A_1,\ldots,A_1,\ldots,A_k,\ldots,A_k$.

For example the mixed volume of $n$ copies of a convex body in $\mathbf{R}^n$ is the usual volume of the body. The mixed volume of $n-1$ copies of a convex body $A$ and one copy of the unit ball is the $n-1$-dimensional volume $V_{n-1}(\partial A)$ of the boundary of $A$ divided by $n$. See for example \cite{burago1988geometric} for introduction and further examples.

For a set $A_1,\ldots,A_k$ of $n\times n$ real symmetric matrices the mixed discriminant $$\det (\underbrace{A_1,\ldots,A_1}_{i_1},\ldots,\underbrace{A_k,\ldots,A_k}_{i_k})$$ is the coefficient $D_{\{i_1,\ldots,i_k\}}$ in the expansion
$$\det( \lambda_1 A_1+\ldots +\lambda_k A_k)=\sum_{\substack{I\subset \mathbf{Z}_{+}^k \\ |I|=n}} {n \choose I} D_I \lambda^I$$

\subsection{Example}

If each $A_i$ is a box $A_i=[0,a_{i1}]\times\ldots\times[0,a_{in}]$ then the mixed volume of the bodies $A_1,\ldots,A_n$ is the permanent of the matrix $(a_{ij})_{i,j=1}^n$ divided by $n!$.

Similarly if $A_i$ is the diagonal matrix with diagonal entries $a_{i1},\ldots,a_{in}$ then the mixed discriminant of the matrices $A_1,\ldots,A_n$ is the permanent of the matrix $(a_{ij})_{i,j=1}^n$  divided by $n!$.

\section{Alexandrov-Fenchel Inequality and Its Analogue in Linear Context}

Alexandrov-Fenchel inequality states that for a set of $n$ bodies $A_1,\ldots,A_n$ in $\mathbf{R}^n$
$$V(A_1,A_2,A_3,\ldots,A_n)^2\geq V(A_1,A_1,A_3,\ldots,A_n)V(A_2,A_2,A_3,\ldots,A_n)$$

This inequality generalizes many known inequalities for convex bodies, like isoperimetric inequality $$\frac{V_n(A)^{1/n}}{V_{n-1}(\partial A)^{1/( n-1)}}\leq\frac{V_n(B^n)^{1/n}}{V_{n-1}(\partial B^n)^{1/(n-1)}}$$  Brunn-Minkowski inequality $$V(A+B)^{1/n}\geq V(A)^{1/n}+V(B)^{1/n}$$ and many others. 

Alexandrov-Fenchel inequality has been first announced by Fenchel \cite{fenchel} and proved by Alexandrov in \cite{alexandrov2} and \cite{alexandrov4}. A simpler proof has been recently found in \cite{mcmullen1993simple} and \cite{timorin1999analogue}.

This inequality is important not only because it is one of the most general inequalities known about convex bodies, but also because of its relations with algebraic geometry. Khovanskii and Tessier have found (see \cite{khovanskiialgebra}, \cite{teissier1979theoreme} and the discussion in \cite{gromov1990convex}) that Alexandrov inequality can be derived from its algebro-geometric analogue $$[D_1,D_2,D_3,\ldots,D_n]^2\geq [D_1,D_1,D_3,\ldots,D_n][D_2,D_2,D_3,\ldots,D_n]$$ where $D_i$ are ample divisors in a smooth irreducible algebraic variety and $[-,\ldots,-]$ is the intersection index of divisors. This algebro-geometric analogue can be derived as a consequence of Hodge index theorem on a smooth irreducible algebraic surface.

Alexandrov-Fenchel inequality has an analogue for  real positive-definite symmetric matrices $A_1,\ldots,A_n$:
$$\det(A_1,A_2,A_3,\ldots,A_n)^2\geq \det(A_1,A_1,A_3,\ldots,A_n)\det(A_2,A_2,A_3,\ldots,A_n)$$

This inequality was first proved in \cite{alexandrov4} and subsequently given a simpler proof in \cite{khovanskii1984analogues}.

Even the most simple cases of these inequalities are extremely useful. For instance when the bodies are boxes with parallel sides (or when the symmetric matrices are diagonal), the inequality on permanents implied by Alexandrov-Fenchel inequality has been used by Falikman \cite{falikman1981proof} in 1979 and by Egorychev in 1980 \cite{egorychev1981solution} to prove a conjecture of van der Waerden that has been open since 1926. Namely he proved that the minimal value of permanent of a doubly stochastic $n\times n$ matrix is attained on the matrix all of whose values are equal to $1/n$.

\section{A Negative Answer to Gromov's Question}

Alexandrov-Fenchel inequalities can be equivalently formulated in terms of the function $f(I)=\log(V_I)$ (or $f(I)=\log(D_I)$) on the discrete simplex $\{I\subset \mathbf{Z}_+^n, |I|=n\}$, where $V_I$ (or $D_I$) are the mixed volumes (or mixed discriminants) appearing in the definition of the mixed volume (or discriminant) above. Namely, assuming in addition that $\log 0= -\infty$, Alexandrov-Fenchel inequality says that the function $f$ is concave on any segment in the discrete simplex that is parallel to one of the sides.

In \cite{gromov1990convex} Gromov asked whether it was true that the function $f$ is concave on the discrete simplex.

In case $n=3$ this generalization amounts to the inequality
$$V(A_1,A_2,A_3)^3\geq V(A_1,A_1,A_2)V(A_2,A_2,A_3)V(A_3,A_3,A_1)$$

 We claim that this inequality fails even when the bodies are boxes with sides parallel to the axes.

Namely we can take $A_1=[0,1]\times[0,1]\times\{0\}$, $A_2=[0,1]\times\{0\}\times[0,5]$ and $A_3=\{0\}\times[0,1/3]\times[0,1]$.

Then $V(A_1,A_2,A_3)=4/9$, $V(A_1,A_1,A_2)=5/3$, $V(A_2,A_2,A_3)=5/9$, $V(A_3,A_3,A_1)=1/9$ and $(4/9)^3< 5/3\cdot 5/9 \cdot 1/9$.

\section{Acknowledgements}

I would like to thank Askold Khovanskii for introducing me to the theory of mixed volumes and Dmitry Faifman and Yevgeny Liokumovich for inspiring discussions.

\bibliography{alexandrov_fenchel}

\begin{thebibliography}{McM93}

\bibitem[Ale37]{alexandrov2}
A.~D. Alexandrov.
\newblock Zur {Theorie} der {Gemischten} {Volumina} von konvexen {K\"{o}rpern}
  {II}; neue ungleichungen zwischen den gemischten {Volumina} und ihren
  {Anwendungen}.
\newblock {\em Math. Sbomik, NS}, 2:1205--1238, 1937.

\bibitem[Ale38]{alexandrov4}
A.~D. Alexandrov.
\newblock Zur {Theorie} der {Gemischten} {Volumina} von konvexen {K\"{o}rpern}
  {IV}; die gemischten {Diskriminanten} und die gemischten {Volumina}.
\newblock {\em Math. Sbomik, NS}, 3:227–251, 1938.

\bibitem[BZ88]{burago1988geometric}
Y.D. Burago and V.A. Zalgaller.
\newblock {\em Geometric Inequalities}.
\newblock Grundlehren der mathematischen Wissenschaften Series. Springer, 1988.

\bibitem[Ego81]{egorychev1981solution}
GP~Egorychev.
\newblock The solution of van der {Waerden}'s problem for permanents.
\newblock {\em Advances in math}, 42:299--305, 1981.

\bibitem[Fal81]{falikman1981proof}
D.I. Falikman.
\newblock Proof of the van der {Waerden} conjecture regarding the permanent of
  a doubly stochastic matrix.
\newblock {\em Mathematical Notes}, 29(6):475--479, 1981.

\bibitem[Fen36]{fenchel}
W.~Fenchel.
\newblock G\'{e}n\'{e}ralisations du th\'{e}or\`{e}me de brunn et minkowski
  concernant les corps convexes.
\newblock {\em C. R. Acad. Sci. Paris}, 203:764--766, 1936.

\bibitem[Gro90]{gromov1990convex}
M.~Gromov.
\newblock Convex sets and {K\"{a}hler} manifolds.
\newblock {\em Advances in Differential Geometry and Topology, ed. F. Tricerri,
  World Scientific, Singapore}, pages 1--38, 1990.

\bibitem[Kho]{khovanskiialgebra}
A.G. Khovanskii.
\newblock Algebra and mixed volumes. {Appendix} 3 in: [{BZ88}].

\bibitem[Kho84]{khovanskii1984analogues}
A.G. Khovanskii.
\newblock Analogues of the {Aleksandrov-Fenchel} inequalities for hyperbolic
  forms.
\newblock In {\em Soviet Math. Dokl}, volume~29, pages 710--713, 1984.

\bibitem[McM93]{mcmullen1993simple}
P.~McMullen.
\newblock On simple polytopes.
\newblock {\em Inventiones mathematicae}, 113(1):419--444, 1993.

\bibitem[Tei79]{teissier1979theoreme}
B.~Teissier.
\newblock Du th{\'e}oreme de l'index de {Hodge} aux in{\'e}galit{\'e}s
  isop{\'e}rim{\'e}triques.
\newblock {\em CR Acad. Sci. Paris Ser. AB}, 288(4), 1979.

\bibitem[Tim99]{timorin1999analogue}
V.A. Timorin.
\newblock An analogue of the {Hodge-Riemann} relations for simple convex
  polytopes.
\newblock {\em Russian Mathematical Surveys}, 54:381, 1999.

\end{thebibliography}

\bibliographystyle{alpha}

Yuri Burda

University of Toronto

yburda@math.toronto.edu
\end{document}